\pdfoutput=1  
\documentclass[11pt,a4paper]{amsart}
\usepackage[T1]{fontenc}
\usepackage[utf8]{inputenc}
\usepackage{lmodern}    
\usepackage{amsmath,amssymb,amsthm}
\usepackage{microtype}
\usepackage[colorlinks=true,linkcolor=blue!50!black,citecolor=blue!50!black,urlcolor=blue!50!black]{hyperref}
\usepackage{xcolor}
\usepackage{booktabs}
\usepackage{array}
\usepackage{alltt}      
\usepackage{stmaryrd}   
\usepackage{tikz}       
\usetikzlibrary{arrows.meta,positioning,decorations.pathreplacing}

\newcommand{\isa}[1]{\texttt{#1}}
\newcommand{\artifactAurl}{https://isa-afp.org/entries/HOL_in_HOL_Deep.html}
\newcommand{\artifactBurl}{https://isa-afp.org/entries/HOL_in_HOL_Deep_Infinity_ZFC.html}
\newcommand{\artA}[2]{\href{\artifactAurl}{\normalfont(Isabelle: \isa{#1}, theory \isa{#2} of \ensuremath{\mathcal{A}_1})}}
\newcommand{\artB}[2]{\href{\artifactBurl}{\normalfont(Isabelle: \isa{#1}, theory \isa{#2} of \ensuremath{\mathcal{A}_2})}}

\newcommand{\Ineq}{\mathrm{Ineq}}
\newcommand{\Exn}[1]{\exists_{\ge #1}}
\newcommand{\DInf}{\mathrm{DInf}}
\newcommand{\ACrel}{\mathrm{AC}^{\mathrm{rel}}}
\newcommand{\NK}{\mathrm{NK}}

\newcommand{\der}{\vdash}
\newcommand{\fder}{\Vdash}          
\newcommand{\sem}[1]{\models_{#1}}
\newcommand{\Mbfb}{\mathcal{M}_{\beta\mathfrak{f}\mathfrak{b}}}
\newcommand{\D}{\mathcal{D}}
\newcommand{\ii}{\iota}
\newcommand{\oo}{o}
\newcommand{\Vs}{V}

\theoremstyle{plain}
\newtheorem{theorem}{Theorem}[section]
\newtheorem{proposition}[theorem]{Proposition}

\newtheorem{corollary}[theorem]{Corollary}
\theoremstyle{definition}

\theoremstyle{remark}
\newtheorem{remark}[theorem]{Remark}
\newtheorem{convention}[theorem]{Convention}

\begin{document}

\title[HOL in HOL: soundness, completeness and consistency]{HOL in
  HOL: mechanized soundness, completeness and consistency, with and
  without infinity}

\author{Christoph Benzm\"uller}
\address{Otto-Friedrich-Universit\"at Bamberg and Freie Universit\"at Berlin, Germany}
\email{christoph.benzmueller@uni-bamberg.de}

\author{Daniel Kirchner}
\address{Otto-Friedrich-Universit\"at Bamberg, Germany}
\email{daniel.kirchner@uni-bamberg.de}


\begin{abstract}\vspace*{-2pt}
What can a logic certify about a copy of itself?  For Church's simple
theory of types, formalized as a deep embedding in Isabelle/HOL with
Henkin semantics in combination with a natural-deduction calculus, the
answer has two
tiers.  On the first tier, plain Isabelle/HOL certifies the logic without an
axiom of infinity, by a standard model with one individual.  On the
second tier, the logic with a Dedekind-style axiom of infinity is
certified only in Isabelle/HOL extended by a large enough
set-theoretic universe.  Both tiers are addressed in Harrison's
self-verification of HOL Light.
What the present formalization adds is Henkin completeness, and with
it three things: a semantics that matches the calculus exactly,
compactness and model existence on the first tier, and the boundary
between the tiers.  Comparing two ways of adding infinity to the object logic shows
where the boundary lies.  One is an infinite domain
of individuals, expressed by the infinitely
many first-order sentences ``there are at least $n$ individuals'',
each of which has a finite model; the other is the Dedekind axiom, a
single sentence with no finite model, and there the first tier
ends.  In the calculus the axiom is strictly stronger than the first-order
sentences together: it derives each of them, but they do not derive
it.  The witness is a countable Henkin model in which
every sentence holds and the axiom fails, because its function domain
is too sparse for an injective, non-surjective self-map.  On the second tier, the universe serves one purpose only: to prove
the Dedekind axiom consistent.  Everything that follows from that consistency is obtained in plain
Isabelle/HOL: completeness, and even the countable Henkin model of
the axiom.  This rests on a Harrison-style universe axiom, realized in
Paulson's \isa{ZFC\_in\_HOL} with the domain of every type in
$\Vs_{\omega+\omega}$; with it, full classical higher-order logic with
infinity and relational choice is consistent relative to set theory.
\end{abstract}

\maketitle

\section{Introduction}\label{sec:intro}

A theory that interprets arithmetic cannot prove its own consistency;
this is G\"odel's second incompleteness theorem.  Yet proof assistants
are routinely trusted with the verification of their own logical cores.
Harrison proved the kernel of HOL Light sound and consistent twice:
for the kernel without its axiom of infinity inside full HOL Light,
and for the kernel with it inside HOL Light extended by one axiom
postulating a large universe~\cite{Harrison06}.  Kumar, Arthan, Myreen
and Owens, and later the Candle project, extended the soundness result to HOL's
definitional mechanisms, with the set-theoretic universe as an
explicit assumption of the soundness theorem rather than an
axiom~\cite{Kumar16,Candle}.  Between G\"odel's limit and these results lies a
question that this paper answers by machine-checked theorems, for
Church's simple theory of types~\cite{Church40}: what exactly can a
copy of the logic certify about the logic, and at what point does the
certifying side have to become stronger than the certified one?
Harrison's two results fix the two tiers of the answer for standard
semantics and soundness.  This paper transfers his division to a deep embedding with Henkin
semantics, where a completeness theorem is available, and
completeness has three consequences here.  First, the calculus
acquires a semantics that matches it exactly: $\NK$ proves exactly
what holds in all Henkin models, so every statement about derivability
below is also a statement about all Henkin models.  Second,
compactness and model existence belong to the first tier: both are
theorems of plain Isabelle/HOL, and even the countable Henkin model of
the Dedekind axiom of infinity is built there, from its consistency alone.
Third, both sides of the boundary become visible: finite models show where the
first tier ends, and the Henkin model built by the completeness proof
shows that the two familiar ways of adding infinity, an infinite
domain of individuals and a Dedekind-style axiom, differ there.

The setting is a \emph{deep embedding}: Church's type theory, the
logic underlying Isabelle/HOL~\cite{Isabelle}, HOL4~\cite{HOL4}, HOL
Light~\cite{HOLLight} and the automated provers of the TPTP library's higher-order
format THF~\cite{THF,BenzmullerMiller14}, is the \emph{object} logic, formalized
as mathematical data inside Isabelle/HOL, the \emph{meta}-logic.  Its
syntax is represented by a datatype, and the Henkin-style
semantics~\cite{Henkin50,BBK04} and the natural-deduction calculus $\NK$ of
Benzm\"uller, Brown and Kohlhase (BKK)~\cite{BBK04} are defined over that datatype; object and
meta-logic are kept strictly apart.  What is formalized is the
monomorphic core common to these systems: primitive equality and typed description, no
polymorphism, no choice operator, and no axiom of infinity unless one
is added explicitly.  The calculus is proved sound and Henkin-complete
(Section~\ref{sec:embedding}) and consistent
(Section~\ref{sec:tier-one}).  The two tiers are distinguished by the meta-theory required.  A
statement is on the first tier if plain Isabelle/HOL proves it, be it
a consistency, completeness or model existence statement, and on the
second tier if it is proved only with a large enough set-theoretic
universe as an additional type (Sections~\ref{sec:universe} and~\ref{sec:harrison}).

\subsection*{Tier one: consistency without infinity}  Without an
object-level axiom of infinity, $\NK$ is consistent by a standard
model with finite domains: a single individual suffices, the
signature may be finite, and the model is defined without a choice
operator
(Section~\ref{sec:tier-one}).  No set theory is involved, and G\"odel's second incompleteness
theorem does not apply, since the object logic interprets no
arithmetic; the meta-level infinity of Isabelle/HOL enters only
through the syntax and the datatype inside which the finite models
are built (Section~\ref{sec:tier-one}).  For many uses of higher-order logic that assume no object-level axiom
of infinity, as in the LogiKEy methodology of semantic
embeddings~\cite{LogiKEy}, this first tier already certifies the logic
in use.

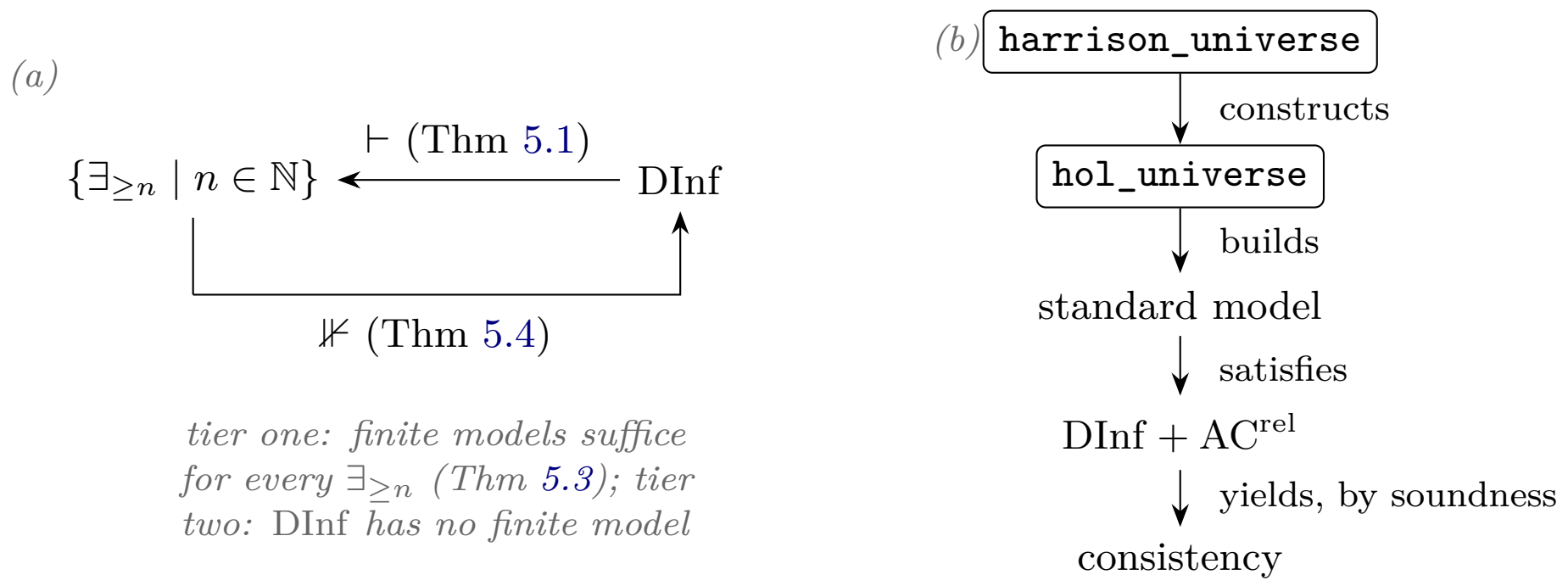
\begin{figure}[t]
\centering
\begin{tikzpicture}[baseline=(current bounding box.center),
  >={Stealth[length=1.8mm]}, font=\footnotesize,
  part/.style={font=\scriptsize\itshape, text=black!60}]
  \node[part] at (0.35,3.10) {(a)};
  \node (exn)  at (1.6,2.30) {$\{\Exn{n} \mid n \in \mathbb{N}\}$};
  \node (dinf) at (5.4,2.30) {$\DInf$};
  \draw[->] (dinf) -- node[midway, above]{$\vdash$ (Thm~\ref{thm:dinfdistinct})} (exn);
  \draw (exn.south) -- ++(0,-0.6) coordinate (a1);
  \draw[->] (a1) -- node[midway, below=1pt]{$\nVdash$ (Thm~\ref{thm:distinctnotdinf})}
    (a1 -| dinf.south) -- (dinf.south);
  \node[part, text width=6.0cm, align=center] at (3.5,-0.05)
    {tier one: finite models suffice for every $\Exn{n}$
     (Thm~\ref{thm:distinctcon}); tier two: $\DInf$ has no finite model};
\end{tikzpicture}%
\hspace{0.6cm}%
\begin{tikzpicture}[baseline=(current bounding box.center),
  >={Stealth[length=1.8mm]}, font=\footnotesize,
  iface/.style={draw, rounded corners=2pt, inner sep=3.5pt,
                font=\footnotesize\ttfamily},
  part/.style={font=\scriptsize\itshape, text=black!60},
  lbl/.style={font=\scriptsize, anchor=west}]
  \node[part] at (-1.75,4.35) {(b)};
  \node[iface] (h)  at (0,4.35) {harrison\_universe};
  \node[iface] (u)  at (0,3.30) {hol\_universe};
  \node (s)   at (0,2.30) {standard model};
  \node (ax)  at (0,1.30) {$\DInf + \ACrel$};
  \node (con) at (0,0.30) {consistency};
  \draw[->] (h) -- (u);
  \draw[->] (u) -- (s);
  \draw[->] (s) -- (ax);
  \draw[->] (ax) -- (con);
  \node[lbl] at (0.17,3.83) {constructs};
  \node[lbl] at (0.17,2.80) {builds};
  \node[lbl] at (0.17,1.80) {satisfies};
  \node[lbl] at (0.17,0.80) {yields, by soundness};
\end{tikzpicture}
\caption{(a)~The two senses of infinity of
  Section~\ref{sec:infinity}: the first-order sentences $\Exn{n}$
  (``at least $n$ individuals'') and the Dedekind axiom $\DInf$.
  Arrows denote derivability in $\NK$ (Convention~\ref{conv:relations});
  the lower arrow is a non-derivability result, witnessed by a Henkin
  model of all the sentences (Theorem~\ref{thm:henkinsep}).  (b)~The
  construction chain behind the tier-two consistency theorem
  (Theorems~\ref{thm:harrison}, \ref{thm:vmodel}
  and~\ref{thm:confullrel}).  $\ACrel$ is the relational choice scheme
  of Section~\ref{sec:universe}; the boxed nodes are the two interfaces
  of $\mathcal{A}_2$, Isabelle locales, that is, named bundles of
  parameters and assumptions; the universe $\Vs$ appears only as one
  instance of the topmost one.}
\label{fig:chain}
\end{figure}

\subsection*{Tier two: consistency with infinity}  A Dedekind-style axiom of infinity $\DInf$ (some self-map of the
individuals is injective but not surjective) has no finite models, so the first tier
cannot certify its consistency.  Section~\ref{sec:boundary} argues that plain Isabelle/HOL cannot prove
the extension by infinity and choice consistent.  The two ways it
offers fail: a standard model needs a carrier type of cardinality at
least $\beth_\omega$, and plain Isabelle/HOL proves no type that
large; a Henkin model is built by the completeness proof only from
the very consistency to be shown.  By G\"odel's second incompleteness
theorem no other way succeeds, since $\NK$ with infinity and choice
interprets the monomorphic core of the meta-logic (this is a
classical argument, not formalized here).  So the meta-theory has to be strengthened, and the strengthening we
need is small: a Harrison-style universe axiom, formalized as the
locale \isa{harrison\_universe}.  From it we construct an abstract
universe interface, the locale \isa{hol\_universe}, and we realize
both in the set-theoretic universe $\Vs$ of Paulson's
\isa{ZFC\_in\_HOL}~\cite{ZFCinHOL}.  The result is a machine-checked
consistency theorem for \emph{full} classical higher-order logic in
Church's sense, extensionality, typed description, the Dedekind axiom
and a relational axiom of choice, relative to set theory
(Sections~\ref{sec:universe} and~\ref{sec:harrison}; the construction
chain is Figure~\ref{fig:chain}(b)).  Choice adds no further
requirement: the relational choice scheme holds in every standard model,
whose full function spaces contain the choice function that the
meta-logic's own choice provides (Section~\ref{sec:choice}).  A rank analysis shows the domain of every type in the concrete model
to be an element of $\Vs_{\omega+\omega}$, that is, to lie within
finitely many power-set steps above the natural numbers
(Section~\ref{sec:rank}).  This is the part of the set-theoretic
universe that Zermelo set theory, ZF without replacement, already
provides; that the construction needs no replacement is an
interpretation of the rank bound, not a formalized reduction to
Zermelo set theory.

\subsection*{The boundary, and what stays below it}  Where exactly
does the first tier end?  Section~\ref{sec:infinity} compares the two
senses of adding infinity of Figure~\ref{fig:chain}(a): infinity as a
semantic property of the individual domain, axiomatized by the first-order sentences ``there are at least $n$
individuals'', and infinity as the single syntactic axiom $\DInf$.
The sentences are on the first tier: $\NK$ is consistent with all of
them by compactness over arbitrarily large finite models, without
ever constructing an infinite model.  $\DInf$ derives every one of
them, but they do not derive $\DInf$: a Henkin model of all of them
refutes it, because a sparse function domain need not contain the
witnessing self-map.  In Henkin semantics the Dedekind axiom is therefore
strictly stronger than an infinite individual domain, while at
standard models the two coincide.  The first tier holds more than the consistency of the first-order
sentences.  Completeness of
$\NK$ with $\DInf$ as a hypothesis is a theorem of plain
Isabelle/HOL, and so is model existence: every consistent sentence
over a countable signature, an axiom of infinity included, has a
Henkin model with countable total domain
(Section~\ref{sec:completeness-vs-standard}).  Only the consistency of the extension requires the stronger
meta-theory, and nothing else is taken from it: the countable
Henkin model of $\DInf$ is built in plain Isabelle/HOL from that
consistency alone.

\subsection*{Contributions}
The contribution of this paper is a precise answer to the opening
question: what higher-order logic can certify about a copy of itself,
where the boundary between the two tiers runs, and how little has to
be added to reach the second tier.  Every theorem, proposition and corollary of this paper is a
machine-checked theorem, and the machine-checking has made the answer
more precise: every hypothesis is explicit, and several of the
distinctions below (the two senses of infinity, the per-type rank
bound, the exact role of choice) emerged where a plausible stronger
claim could not be proved.  The theorems are contained in two entries of the
Archive of Formal Proofs (AFP): \isa{HOL\_in\_HOL\_Deep}~\cite{HOLinHOLDeep},
extended for this paper and called $\mathcal{A}_1$ below, and its
companion \isa{HOL\_in\_HOL\_Deep\_Infinity\_ZFC}~\cite{InfinityZFC},
called $\mathcal{A}_2$.  The metatheoretic steps that are not
formalized (Section~\ref{sec:boundary}, Table~\ref{tab:matrix},
Remark~\ref{rem:foundations}) are marked as such wherever they occur.
In detail, we provide:
\begin{enumerate}
\item consistency of $\NK$ without infinity, by a finite standard
  model, with a footprint audit: no object-level infinity, no set
  theory, no choice operator in the definition of the model, and a possibly finite
  signature (Section~\ref{sec:tier-one});
\item the plain-HOL results: soundness, Henkin completeness relative
  to hypotheses at every infinite signature, semantic compactness, and
  model existence with a countable witness
  (Sections~\ref{sec:embedding} and~\ref{sec:completeness-vs-standard});
\item the boundary, measured by two senses of infinity: consistency
  of the first-order sentences ``there are at least $n$ individuals''
  by compactness over finite models; their derivability from $\DInf$;
  non-derivability of $\DInf$ from them, with a Henkin countermodel
  satisfying all of them; and the characterization of satisfaction of
  $\DInf$ in any Henkin model (Section~\ref{sec:infinity});
\item an abstract universe interface with a carrier-generic
  consistency theorem for the extension by infinity and choice
  (Section~\ref{sec:universe});
\item a Harrison-style interface from which the universe interface is
  constructed, and its realization in \isa{ZFC\_in\_HOL}, yielding
  the consistency of full classical HOL in Church's sense
  (Sections~\ref{sec:universe} and~\ref{sec:harrison});
\item a type-by-type rank bound: every domain of the concrete model is
  an element of $\Vs_{\omega+\omega}$ (Section~\ref{sec:rank}).
\end{enumerate}
Items (1), (4), (5) refine Harrison's two tiers; (2), (3), (6) have
no counterpart there.

\section{The underlying deep embedding}\label{sec:embedding}

We recall from $\mathcal{A}_1$~\cite{HOLinHOLDeep} what the later
sections use; for the full development see the entry itself.  The object logic has base types $\ii$ (individuals) and
$\oo$ (truth values) and function types $\sigma \Rightarrow \tau$.

The syntax, after BKK~\cite{BBK04}, is formalized in a \emph{locally
nameless} representation~\cite{Chargueraud12}: a bound variable is
not a name but a number, counting the binders between its occurrence
and the binder that captures it, while a free variable is a named
atom; terms that differ only in the names of bound variables are
therefore identical.  The signature is a set of typed constants
(\emph{parameters}) drawn from a meta-level type $p$;\footnote{In the
  artifact $p$ is a type variable (\isa{'p}), where required with the
  class constraints \isa{infinite} or \isa{countable}; every
  statement is thus polymorphic in the signature.} BKK's
identification of $\alpha$-equivalent terms thereby becomes a literal
equality.
Concretely, a term is a value of the datatype \isa{tm}, with
constructors for bound variables (de Bruijn indices), free variables
(a name and a type), parameters (a signature constant and a type),
the logical constants $\neg$, $\vee$, $\Pi_\sigma$, $\iota_\sigma$
and $=_\sigma$, application, and abstraction with an annotated domain
type.  Well-formedness at a type $\sigma$ is an inductive predicate.
The binder notation of the displayed formulas, $\forall X_\sigma.\,A$
and $\exists X_\sigma.\,A$, is defined in the usual way by abstracting a named free
variable, so that every displayed formula is literally a value of
\isa{tm} (Remark~\ref{rem:reading}).  The signature includes primitive equality and, following
Andrews~\cite{Andrews72} and going beyond BKK, typed description
operators.  The calculus is BKK's $\NK_{\beta\mathrm{fb}}$ (their
Definition~7.1 and Figures~6, 7 and~9), fifteen rules: hypothesis,
$\beta$-conversion, introduction and elimination for $\neg$, $\vee$
and $\Pi$, classical contradiction, functional and Boolean
extensionality, two rules for primitive equality, and, after Andrews
but at every type, one rule for description.  The semantics is rendered abstractly, for the class
$\Mbfb$ of $\Sigma$-Henkin models (``$\Sigma$'' is BKK's name for the
signature, our $p$) over an arbitrary \emph{value carrier} $u$, the
meta-level type inside which all domains live.  A model is a frame of
typed domains $\D_\sigma$ with application, in which every well-formed
term denotes and the logical constants behave as expected;
\emph{standard} models are the special case in which every function
domain $\D_{\sigma\Rightarrow\tau}$ is the full function space.

The results of this paper are stated with three consequence
relations, two syntactic and one semantic, fixed in the following
convention.

\begin{convention}\label{conv:relations}
For a set $\Phi$ of object formulas and a formula $A$:
\begin{itemize}
\item $\Phi \der A$ denotes derivability in $\NK$ from the \emph{whole}
  context $\Phi$ \artA{bprov}{Calculus};
\item $\Phi \fder A$ denotes the \emph{finitary} consequence in the style
  of Andrews: some finite $\Phi_0 \subseteq \Phi$ satisfies
  $\Phi_0 \der A$ \artA{fprov}{Calculus};
\item $\Phi \sem{u} A$ denotes semantic consequence over all models of
  $\Mbfb$ with value carrier $u$ \artA{bkk\_consequence}{Semantics}.
\end{itemize}
$\mathrm{con}\,\Phi$ abbreviates $\Phi \nvdash \bot$.  Over infinite
signatures, $\fder$ and $\der$ provably coincide on finite contexts
and on \emph{pure} contexts, those whose formulas contain no
parameters \artA{fprov\_eq\_bprov}{Calculus}.
\end{convention}

\begin{remark}[How to read the formal statements]\label{rem:reading}
Every displayed object formula is, in the artifact, literally a value of
an Isabelle datatype: the type \isa{tm} of locally nameless terms over
the signature type $p$ (theory \isa{Syntax}).  The axiom $\DInf$ of Section~\ref{sec:boundary}, for instance, is
the Isabelle definition \isa{DInf}, a term of type \isa{tm} written
with the binder notation in the form displayed there.  That this
term is a closed, well-formed formula is not stipulated, as it would
be on paper, but proved (\isa{cwff\_DInf}, theory \isa{NK\_Infinity}).  On the semantic side, ``$A$ is satisfied in $M$ under assignment
$\xi$'' means that the denotation of $A$ under $\xi$ in $M$ is the
distinguished truth value, in the artifact the equation
$\isa{den}\,A\,\xi = \isa{Tv}$.  An assignment is always
type-respecting: each variable of type $\sigma$ receives a value in
$\D_\sigma$ (\isa{asg}; both in theory \isa{Semantics}).  Each
principal statement below is a natural-language reading of a formal
statement in the artifact, with every hypothesis shown, and carries a
pointer of the form ``(Isabelle: \emph{name}, theory \emph{T} of
$\mathcal{A}_i$)'' to it.
\end{remark}

The main theorem of $\mathcal{A}_1$ states soundness and completeness
under two hypotheses, an infinite signature and a value carrier into
which the signature embeds; every later section uses it in this form.

\begin{theorem}[Soundness and completeness, hypothesis-relative]
\label{thm:soundcomplete}
Let $\Phi$ be any set of well-formed object formulas over an
\emph{infinite} signature $p$
and $A$ well-formed, and let the value carrier $u$ admit an injection
$p \hookrightarrow u$.  Then
\[
  \Phi \fder A \quad\Longleftrightarrow\quad \Phi \sem{u} A .
\]
No purity, finiteness, countability, or closedness condition is
required of $\Phi$ or $A$.
\artA{sound\_and\_complete\_hyps}{Main\_Results}
\end{theorem}

\begin{corollary}[Semantic compactness]\label{cor:compact}
Let $\Phi$, $A$ and $u$ be as in Theorem~\ref{thm:soundcomplete}.  If
$\Phi \sem{u} A$ then, for every value carrier $v$, some finite
$\Phi_0 \subseteq \Phi$ satisfies $\Phi_0 \sem{v} A$.
\artA{NK\_consequence\_compact}{Main\_Results}
\end{corollary}

The proof uses no ultraproducts.  Completeness turns $\Phi \sem{u} A$
into a derivation of $A$ from some finite $\Phi_0 \subseteq \Phi$,
and soundness, which holds over every carrier, turns that derivation
into $\Phi_0 \sem{v} A$ for every $v$.\footnote{For a carrier $v$ too
  small to admit any model of $\Mbfb$, such as \isa{unit},
  $\Phi_0 \sem{v} A$ holds vacuously.}  The same $\Phi_0$ thus serves
all carriers.  The formal corollary cannot say so, because
Isabelle/HOL quantifies over types only at the outermost level: $v$
is a schematic type variable, and no quantifier over it can follow
the existential over $\Phi_0$.  The corollary is therefore one
theorem per carrier $v$; the uniformity is recorded by the syntactic
statement $\Phi \fder A$ of Theorem~\ref{thm:soundcomplete}, whose
witness $\Phi_0$ mentions no carrier.

\section{Tier one: consistency without infinity}\label{sec:tier-one}

The first consistency theorem of the development needs nothing beyond
plain Isabelle/HOL.

\begin{theorem}[Consistency of $\NK$]\label{thm:consistency}
$\nvdash \bot$, at every signature $p$; consequently no formula is
derivable together with its negation.
\artA{consistency}{Main\_Results}, \artA{no\_contradiction}{Main\_Results}
\end{theorem}

The proof exhibits a model and appeals to soundness: whatever is
derivable is true in the model, and $\bot$ is not.  The model is
standard and finite: one individual, two truth values, and every
function domain the full function space, finite because its
argument and result domains are, represented by finite graphs.  This is Andrews' consistency argument for
$\mathcal{Q}_0$ (Theorem~5403 of~\cite{Andrews02}).  The logical
constants are not defined individually: the frame is an instance of
an abstract interface (\isa{lambda\_universe}) over which negation,
disjunction, quantification, equality and description are
\emph{constructed} \artA{concrete\_lambda\_universe}{Consistency}; the
description operator is interpreted by Isabelle's own definite
description \isa{THE}.  The artifact builds
such a model for any number $k > 0$ of individuals; the larger
instances serve the finite-model arguments of
Section~\ref{sec:infinity}, and their parameter interpretation can
distinguish any prescribed finite family of constants.  The
one-element model refutes every sentence without a finite model, such
as an axiom of infinity.  The first tier asserts the consistency of no
such sentence, which is why it needs neither an infinite model nor a
stronger meta-theory.

What the proof does \emph{not} use is recorded by a footprint audit
(theory \isa{Footprint\_Checks}), which instantiates the consistency
theorem at the finite parameter types \isa{bool} and \isa{unit}
\artA{consistency\_bool\_params}{Footprint\_Checks}.  The object logic
may have a finite signature, the model has finite domains at every type,
and no set theory is involved.

No choice operator occurs in the model itself, while the proofs about
it inherit choice as a dependency on the meta-logic.  The defining
equations of the model, for its domains, application,
abstraction, parameter interpretation and assignment, contain neither
\isa{SOME} nor \isa{THE}; the description operator comes from the
interface, where it is Isabelle's definite description \isa{THE}, not
a choice operator.  Theory \isa{Footprint\_Deps} checks the first
claim mechanically.  The same theory audits the proof terms of the
consistency theorem and finds the meta-logic's choice principles in
four places: the datatype package, the lemmas about fresh variable
names in the syntax, the evaluation of abstraction and quantification
in the semantics, and the canonical assignment of the soundness
proof.  The footprint claim is therefore ``no choice in the model'',
not ``no choice in the proof'' (Remark~\ref{rem:foundations}).

The kind of statement obtained deserves to be made precise.  The
meta-logic Isabelle/HOL has an axiom of infinity; the object logic has
none.  The meta-level infinity is used in two places.  The syntax uses
it: the term type and the inductively defined set of derivations are
infinite, as they must be in any meta-theory that speaks about
derivations of arbitrary length.  The models use it too, as formalized: their value carrier $u$ is one
datatype \isa{ival}, a recursive type and hence infinite, inside which
every finite model is built, and the proofs about the models use the
natural numbers freely.  The content of each instance,
by contrast, uses nothing beyond finite objects: a single model with
finite domains, in which evaluation is decidable, and a soundness
proof by induction over derivations.  A finite fragment of \isa{ival}
would suffice for any single derivation, which mentions finitely many
types; this is a metatheoretic
remark (not formalized here): the model for
all types is infinite as a whole, and the formalization does not
isolate such fragments.  So the consistency proof is carried out in a logic with an axiom of
infinity, about a logic without one.  No logic proves its own consistency here, and G\"odel's second
incompleteness theorem does not apply to the object logic: a theory with a standard model all of whose domains are finite,
one individual and full function spaces above it, interprets no
arithmetic.

For many applications this is no restriction.  Where HOL serves as a
representation language---as in the LogiKEy methodology of semantic
embeddings~\cite{LogiKEy}, whose embedded object logics assume no
object-level axiom of infinity---the monomorphic core of the logic
employed (TH0, the monomorphic dialect of the TPTP THF
family~\cite{THF}, without its choice operator) is the one certified
here.  The certificate
makes no claim about Isabelle/HOL's polymorphic definitional
superstructure, whose consistency is established
separately~\cite{KuncarPopescu19}.

\section{The boundary: infinity on the object side}\label{sec:boundary}

The \emph{Dedekind-style axiom of infinity}~\cite{Dedekind1888} is the
single \emph{closed,
pure} object formula
\[
  \DInf \;=\;
  \exists F_{\ii\Rightarrow\ii}.\;
  \bigl(\forall X\, Y.\; F X = F Y \supset X = Y\bigr)
  \;\wedge\;
  \bigl(\exists C.\; \forall Z.\; \neg\, (F Z = C)\bigr)
\]
\artA{DInf}{NK\_Infinity}, asserting the existence of an injective, non-surjective self-map of the
individuals.  It mentions no parameters at all.  Up to notation and elementary
logical equivalence, $\DInf$ is the axiom of infinity of HOL Light
(\isa{INFINITY\_AX}~\cite{Harrison06}): HOL Light and HOL4 postulate
this sentence, and Isabelle/HOL an equivalent skolemized form.

\begin{proposition}[Finite refutation]\label{prop:dinfref}
Every Henkin model with a finite individual domain refutes $\DInf$,
under every type-respecting assignment.
\artA{DInf\_refuted\_finite}{NK\_Infinity}
\end{proposition}

Proposition~\ref{prop:dinfref} is the pigeonhole principle in semantic
form: an injective self-map of a finite set is surjective.  In
particular, the finite models of Section~\ref{sec:tier-one} refute
$\DInf$, so the argument of the first tier does not extend to
$\mathrm{con}\,\{\DInf\}$.

Could an infinite model be built in plain Isabelle/HOL instead?
Plain Isabelle/HOL has infinite types, and every level $\mathbb{N}$,
$\mathbb{N} \Rightarrow \mathrm{bool}$, \dots{} of the full
function-space tower exists as a type.  What it cannot prove is that
\emph{one type carries the whole tower at once}.  The obstruction is
stated in four steps: three observations, on the shape of a
formalized model~(1), on standard models~(2) and on Henkin
models~(3), and one argument that carries the conclusion~(4):

\begin{enumerate}
\item In the deep embedding the object types are values of a
  datatype, not types of the meta-logic.  A formalized model is
  therefore one evaluation function, taking terms of all object types
  to values of one meta-level type $u$, the value carrier, inside
  which all domains live (\isa{ival} for the finite models of
  Section~\ref{sec:tier-one}, Paulson's $\Vs$ for the model of
  Section~\ref{sec:harrison}).  The codomain cannot vary with the
  object type, because the simple types of the meta-logic do not
  depend on values.  The single carrier is thus not part of BKK's
  notion of model, where ``standard'' means only that the function
  domains are full; it is imposed by the formalization in
  simply-typed HOL.
\item Over an infinite individual domain with full function spaces,
  the carrier must have cardinality at least $\beth_\omega$, the
  supremum of the cardinalities of the levels of the tower.  Plain
  Isabelle/HOL cannot prove any of its types to be that large.  For a
  type $\tau$ to carry the tower, each level $\mathbb{N}$,
  $\mathbb{N} \Rightarrow \mathrm{bool}$, \dots{} must inject into
  $\tau$.  In the standard set-theoretic semantics of HOL, with a
  countable type \isa{ind}, every type has cardinality below
  $\beth_\omega$, because each type is built from \isa{ind} by
  finitely many function-space and definitional steps.  So for every
  $\tau$ one of these injections fails in that semantics and is
  therefore, by soundness, not a theorem.
\item Henkin models avoid the cardinality obstruction, but they need
  not satisfy $\DInf$ (Remark~\ref{rem:semantic-gap}).  The canonical
  way of building one that does, the term-model construction of
  Section~\ref{sec:completeness-vs-standard}, takes
  $\mathrm{con}\,\{\DInf\}$ as its \emph{premise}.
\item That no other way succeeds either is G\"odel's second
  incompleteness theorem.  $\NK$ with infinity and choice is Church's
  type theory with infinity and choice, and this theory interprets
  the monomorphic core of the meta-logic: the meta-level type of
  individuals by the object individuals, infinite by $\DInf$, and
  function types by function types.  A consistency proof of $\NK$
  with infinity and choice in plain Isabelle/HOL would therefore
  yield a consistency proof of Isabelle/HOL's own core, which
  G\"odel's theorem excludes.\footnote{The step takes as given that plain
  Isabelle/HOL is conservative over its monomorphic core with choice:
  polymorphism by type instantiation of the proof, the definitional
  mechanisms for constants and types by the result of Kun\v{c}ar and
  Popescu~\cite{KuncarPopescu19}.  We do not examine this here.}  The two theories are in
  fact of the same strength: by the classical analysis both are
  equiconsistent with Mac Lane set theory, a subsystem of Zermelo set
  theory (Mathias~\cite{Mathias01}), a standard account we rely on
  and do not formalize; the consistency of type theory with infinity
  in Zermelo set theory goes back to Kemeny~\cite{Kemeny49}.
  For the extension by infinity alone we make no such claim: the
  argument would need in addition that adjoining the choice scheme
  does not raise consistency strength, provably in the theory itself,
  a classical fact whose proof without replacement is
  delicate~\cite{Mathias01} (Remark~\ref{rem:foundations}).
\end{enumerate}

Steps (1) and (3) appear in the formalization; step (2) is a
cardinality argument about the standard semantics of HOL, carried out
at the meta-level; step (4) is a metatheoretic observation the
artifacts discuss but do not formalize.

The conclusion is that the meta-theory has to be strengthened, and the obstacle
is a familiar one: Harrison's self-verification of HOL Light removes
it by a universe axiom~\cite{Harrison06}; its successors add no axiom
but state soundness under an explicit assumption, a type with a
membership relation satisfying the axioms of set
theory~\cite{Kumar16,Candle}.
Sections~\ref{sec:universe} to~\ref{sec:rank} show how little
strength beyond plain Isabelle/HOL the second tier needs, through two
small interfaces.  In the terms of a hierarchy of infinity axioms, write $\mathrm{HOL}_0$
and $\mathrm{HOL}_1$ for higher-order logic without and with an axiom
of infinity.  The meta-logic Isabelle/HOL is $\mathrm{HOL}_1$: its own
axiom of infinity is $\DInf$ on the type \isa{ind} with the witnesses
named---theory \isa{Nat} postulates constants \isa{Zero\_Rep} and
\isa{Suc\_Rep} with \isa{Suc\_Rep} injective and never
\isa{Zero\_Rep}, the skolemized form mentioned above, interderivable
with $\DInf$ by naming the witnesses and by $\exists$-introduction
(not formalized).  Tier one is $\mathrm{HOL}_1 \vdash
\mathrm{con}(\mathrm{HOL}_0)$, while $\mathrm{HOL}_1 \vdash
\mathrm{con}(\mathrm{HOL}_1)$ is not available.  Tier two, as formalized, is $\mathrm{HOL}_1 + \isa{ZFC\_in\_HOL}
\vdash \mathrm{con}(\mathrm{HOL}_1)$.  Its proof depends on Paulson's
universe only through the seven assumptions of the Harrison-style
interface (Section~\ref{sec:harrison}).  This bounds what the proof
uses, not the consistency strength it needs: that every interpretation
of the interface yields $\mathrm{con}(\mathrm{HOL}_1)$ is a theorem of
the artifact, but the only interpretation it provides is Paulson's
universe $\Vs$, and no weaker one, such as a single strong-limit type
above \isa{ind}, is shown to exist.  The reading ``$\mathrm{HOL}_1$
plus one strong-limit type'', suggested by Harrison's axiom $I$,
therefore describes the interface informally; it is not a theorem of
the artifact.
Sections~\ref{sec:infinity} and~\ref{sec:completeness-vs-standard}
measure where the boundary between the tiers runs;
Sections~\ref{sec:universe} to~\ref{sec:rank} then build the second
tier.

\section{Measuring the boundary: two senses of infinity}\label{sec:infinity}

\subsection{The first-order sentences}
The axiomatization of ``the domain is infinite'' in the language of
pure equality~\cite[III.6.3]{EFT21} consists of the constant-free sentences
\[
  \Exn{n} \;=\; \exists x_1 \dots \exists x_n.\; \textstyle\bigwedge_{i<j}
    \neg\,(x_i = x_j) ,
  \qquad n \in \mathbb{N},
\]
one for each $n$ \artA{ExDistinct}{NK\_Infinity}.  Each has a finite
model, while the set of all of them has none (in general, the class
of infinite structures is $\Delta$-elementary, that of finite
structures is not~\cite[VI.3.4]{EFT21}).  This, and not the order of the language, is
what separates the sentences from $\DInf$.\footnote{A first-order sentence can lie on the same side of the
  boundary as $\DInf$: over one function symbol $f$, ``$f$ is
  injective and not surjective''~\cite[VI.1]{EFT21} is the skolemized
  form of $\DInf$, with the function symbol in place of the
  existentially bound $F$, and it has no finite model either.}  The boundary between the tiers thus runs between sentences
that have finite models and sentences that have none.  Inside $\NK$,
$\DInf$ derives every one of the sentences:

\begin{theorem}[$\DInf$ derives every finite cardinality]\label{thm:dinfdistinct}
Over an infinite signature $p$, for every $n$:\; $\{\DInf\} \der \Exn{n}$.
\artA{axiom\_derives\_\allowbreak every\_\allowbreak finite\_cardinality}{Main\_Results}
\end{theorem}

The derivation is built by meta-level recursion on $n$, one
$\NK$-derivation for each $n$: two parameters serve as eigenvariables
for the self-map $F$ and the point $C$ outside its range, the witnesses
are $C, F\,C, F\,(F\,C), \dots$, and their pairwise distinctness
follows from injectivity.

\subsection{The semantic bridge}
\begin{proposition}[Satisfaction of $\DInf$]\label{prop:dinfsat}
In every Henkin model and under every type-respecting assignment,
$\DInf$ is satisfied
if and only if some element $d$ of the function domain
$\D_{\ii\Rightarrow\ii}$ is injective on $\D_\ii$ and omits a point of
$\D_\ii$.  In a standard model, where that domain is full, this holds
if and only if $\D_\ii$ is Dedekind-infinite as a set, that is,
admits an injective self-map that is not surjective.
\artA{DInf\_sat\_iff}{NK\_Infinity},
\artA{DInf\_sat\_iff\_Dedekind}{NK\_Infinity}
\end{proposition}

\subsection{The first-order sentences are consistent}

\begin{theorem}[Consistency of the first-order sentences]\label{thm:distinctcon}
At every signature $p$, finite or infinite:
$\;\neg\,(\{\Exn{n} \mid n\} \fder \bot)$.
\artA{finite\_cardinalities\_\allowbreak consistent}{Main\_Results}
\end{theorem}

The proof constructs \emph{no} infinite model: a finite part of the
sentences mentions a largest $n$, and the finite standard model of
Section~\ref{sec:tier-one} with at least $n$ individuals satisfies it,
so no finite part derives falsity, which is the statement.  This is
the textbook compactness argument that a theory with arbitrarily large
finite models has an infinite model, carried out with the very same
sentences~\cite[VI.2.2]{EFT21}---except that no infinite model is ever
formed.  The sentences mention no constants, which is why the
signature may be finite; over an infinite signature the context, being
pure, has $\fder$ and $\der$ agree \artA{fprov\_eq\_bprov}{Calculus},
and the theorem may there equally be read as
$\mathrm{con}\,\{\Exn{n} \mid n\}$.

The artifact also carries out the finite-model argument in the
classical form of the compactness argument by expansion with new
constants~\cite{Skolem34,Malcev36,Henkin49,Robinson63}, applied in
type theory in Theorem~5506 of Andrews' \S 55~\cite{Andrews02}.  For an
injective family $f : \mathbb{N} \to p$ of parameters over an infinite
signature, the scheme $\Ineq_f = \{\neg(f_i =_\ii f_j) \mid i \neq j\}$
is consistent by the same finite models
\artA{con\_Ineq\_fprov}{NK\_Infinity}; it derives every $\Exn{n}$
\artA{Ineq\_derives\_distinct\_n}{NK\_Infinity}; and a cut reduces every finitary consequence of the sentences to one of the
scheme \artA{distinct\_scheme\_fprov\_Ineq}{NK\_Infinity}.

\subsection{The first-order sentences do not derive $\DInf$}

\begin{theorem}[The first-order sentences do not derive $\DInf$]\label{thm:distinctnotdinf}
At every signature $p$, finite or infinite:\; $\neg\,(\{\Exn{n} \mid n\} \fder \DInf)$.
\artA{finite\_cardinalities\_\allowbreak do\_not\_\allowbreak derive\_axiom}{Main\_Results}
\end{theorem}

Again a finite-model argument: finitely many of the sentences hold in
a finite standard model, which refutes $\DInf$
(Proposition~\ref{prop:dinfref}), so by soundness they do not derive
it.  With Theorem~\ref{thm:dinfdistinct}, the first-order
axiomatization of ``the domain is infinite'' is thus strictly weaker
in $\NK$ than $\DInf$.  Theorem~\ref{thm:distinctnotdinf} separates
$\DInf$ only from finite parts of the set, and no standard model
separates it from the whole set (Remark~\ref{rem:semantic-gap}); the
separating model is a Henkin model, supplied by the completeness
proof:

\begin{theorem}[Henkin countermodel]\label{thm:henkinsep}
Over every infinite signature $p$ there is a $\Sigma$-Henkin model
whose total domain embeds into the term type, that satisfies every
sentence $\Exn{n}$ and refutes $\DInf$.
\artA{henkin\_scheme\_refutes\_DInf}{NK\_Infinity}
\end{theorem}

The model is the term model that the completeness proof builds to
refute a non-derivable formula; over a countable signature the term
type is countable, and so is the model.  The construction needs a reserve of unused parameters as large as
the signature, yet the theorem carries no such hypothesis, because
the reserve is created inside the signature.  The signature is mapped
injectively into a part of itself whose complement has the
signature's cardinality; the model is built there for the renamed
sentences and pulled back along the renaming, which fixes $\DInf$ and
every $\Exn{n}$, both being pure.
(The formal statement takes an injective family $f$ of parameters,
available over any infinite signature, and records that the model also
satisfies the scheme $\Ineq_f$.)

\begin{remark}[Semantic infinitude does not force $\DInf$]
\label{rem:semantic-gap}
The individual domain of the model of Theorem~\ref{thm:henkinsep} is
infinite, since every $\Exn{n}$ holds in it.  Yet $\DInf$ is \emph{false} there: in a
Henkin model the function domain $\D_{\ii\Rightarrow\ii}$ may be too
sparse to contain any injective non-surjective self-map, and in this one
it is.  Conversely, in any Henkin model satisfying $\DInf$, the witnessing function, sparse domains or not,
acts injectively and non-surjectively on $\D_\ii$, so the individual
domain is Dedekind-infinite
\artA{DInf\_sat\_Dedekind\_infinite}{NK\_Infinity}.  In Henkin
semantics, then, the object-level
axiom is strictly stronger than infinitude of the domain; at
standard models $\DInf$ coincides with Dedekind infinitude of the
domain (Proposition~\ref{prop:dinfsat}).  Table~\ref{tab:matrix} summarizes the
comparison.
\end{remark}

\begin{table}[t]
\caption{Henkin versus standard semantics for the embedded logic.  The
  entries name the results of this paper.  The two failures at standard
  semantics are classical readings (of Propositions~\ref{prop:dinfsat}
  and~\ref{prop:dinfref}, and of G\"odel's first incompleteness theorem),
  not formalized theorems; the compactness failure uses the classical step that an infinite set
  admits an injective non-surjective self-map, which needs choice and
  is not formalized.  $\ACrel$ is the relational choice scheme of
  Section~\ref{sec:universe}.}
\label{tab:matrix}
\small
\begin{tabular}{@{}lp{0.32\linewidth}p{0.35\linewidth}@{}}
\toprule
 & \textbf{Henkin models ($\Mbfb$)} & \textbf{standard models} \\
\midrule
completeness of $\NK$ &
  holds (Theorem~\ref{thm:soundcomplete}) &
  fails classically: standard validity is not recursively
  enumerable (not formalized) \\[2pt]
compactness &
  holds (Corollary~\ref{cor:compact}) &
  fails: $\{\Exn{n} \mid n\} \cup \{\neg\DInf\}$ has no standard model,
  though every finite part does \\[2pt]
$\DInf$ versus infinite $\D_\ii$ &
  strictly stronger: $\DInf$ forces $\D_\ii$ infinite, yet can fail
  with $\D_\ii$ infinite (Theorem~\ref{thm:henkinsep},
  Remark~\ref{rem:semantic-gap}) &
  equivalent to Dedekind infinitude of $\D_\ii$
  (Proposition~\ref{prop:dinfsat}) \\[2pt]
$\ACrel$ &
  not guaranteed: a sparse domain may omit the choice function &
  holds always (Proposition~\ref{prop:acrelsat}) \\
\bottomrule
\end{tabular}
\end{table}
\section{Completeness and model existence remain in plain Isabelle/HOL}
\label{sec:completeness-vs-standard}

Completeness for the extension by $\DInf$ needs nothing beyond plain
Isabelle/HOL, neither the universe interfaces of
Sections~\ref{sec:universe} and~\ref{sec:harrison} nor their
realization in Paulson's set-theoretic universe:

\begin{theorem}[Completeness of the $\DInf$ extension, in plain
    Isabelle/HOL]
\label{thm:completeness-dinf}
For every infinite signature $p$ and every infinite value carrier $u$: if
$A$ is well-formed and $\{\DInf\} \sem{u} A$, then $\{\DInf\} \der A$.
\artB{NK\_completeness\_DInf}{HOL\_Universe}
\end{theorem}

The theorem is an instance of the hypothesis-relative completeness of
$\mathcal{A}_1$ (Theorem~\ref{thm:soundcomplete}; the context
$\{\DInf\}$ is finite and pure, and for finite contexts the injection
hypothesis of that theorem can be dropped---infinitude of signature and
carrier suffices), and the theory proving it imports no set
theory.\footnote{At the level of theories, \isa{HOL\_Universe} imports only
theories of $\mathcal{A}_1$.  The session $\mathcal{A}_2$ as a whole
depends on \isa{ZFC\_in\_HOL}, because its second theory instantiates
the interfaces in $\Vs$; the completeness theorem is independent of
that instantiation.}

Completeness says nothing about the consistency of $\DInf$: were
$\DInf$ inconsistent, both sides of the theorem would hold vacuously,
since a hypothesis without models has every formula as a semantic
consequence, and an inconsistent hypothesis derives every formula.

\emph{Model existence} for the extension is likewise available in
plain Isabelle/HOL, conditional on its consistency.
Over a countably infinite signature, $\mathrm{con}\,\{\DInf\}$ as an
explicit hypothesis yields a countable Henkin model of $\DInf$
\artB{countable\_henkin\_model\_DInf\_of\_con}{HOL\_Universe}; what
plain Isabelle/HOL does not supply is the hypothesis itself
(Section~\ref{sec:boundary}), which this development obtains only in
the stronger meta-theory.  Once the hypothesis is discharged there
(Section~\ref{sec:vmodel}), the model is obtained outright, still
constructed in plain Isabelle/HOL
\artB{countable\_henkin\_model\_DInf}{HOL\_in\_HOL\_Deep\_Infinity\_ZFC}.

\section{Standard models: consistency and choice}\label{sec:universe}

\subsection{Choice in full standard semantics}\label{sec:choice}

Let $\ACrel$ denote the relational choice scheme: for each pair of
types $\sigma, \tau$, the closed object formula
\[
  \ACrel_{\sigma\tau} \;=\;
  \forall R_{\sigma\Rightarrow\tau\Rightarrow\oo}.\;
  \bigl(\forall X_\sigma.\; \exists Y_\tau.\; R\,X\,Y\bigr)
  \;\supset\;
  \bigl(\exists F_{\sigma\Rightarrow\tau}.\; \forall X_\sigma.\;
  R\,X\,(F\,X)\bigr)
\]
\artB{ACrel}{HOL\_Universe}, asserting that every left-total relation,
one relating each $X$ to at least one $Y$, has a choice function $F$.

\begin{proposition}[Choice holds in every standard model]
\label{prop:acrelsat}
Every instance $\ACrel_{\sigma\tau}$ of the relational choice scheme
is satisfied in every standard model, under every type-respecting
assignment.
\artB{ACrel\_sat}{HOL\_Universe}
\end{proposition}

The witnessing function is the semantic $\lambda$-abstraction of a
meta-level Hilbert choice of related targets, which lies in the
function domain because that domain is \emph{full}.  This is a theorem
about \emph{standard} models: in Henkin models, whose sparse function
domains need not contain the choice function, the scheme need not
hold (Table~\ref{tab:matrix}).

\subsection{Consistency from an abstract universe}

The consistency argument of the second tier is stated first for an
arbitrary value carrier, with no set theory involved; a concrete
carrier is constructed only in Section~\ref{sec:vmodel}.

\begin{theorem}[Abstract consistency of full HOL]\label{thm:confullrel}
Suppose some standard model, over an arbitrary value carrier, has a
function domain $\D_{\ii\Rightarrow\ii}$ containing an element that is
injective and not surjective on the individual domain $\D_\ii$.  Then
\[
  \mathrm{con}\,\bigl(\{\DInf\} \cup
    \{\ACrel_{\sigma\tau} \mid \sigma, \tau \text{ types}\}\bigr).
\]
\artB{con\_full\_HOL\_rel}{HOL\_Universe}
\end{theorem}

The proof combines the satisfaction of $\DInf$ in such a model
(Proposition~\ref{prop:dinfsat}), the satisfaction of the choice
scheme (Proposition~\ref{prop:acrelsat}) and soundness.  Since $\NK$
already carries Boolean and functional extensionality and typed
description, the theory whose consistency is certified is full
classical higher-order logic in Church's
sense~\cite{Church40,SEPChurchTT}: $\NK$ extended by Dedekind
infinity and relational choice.  Here ``full'' means that all of
Church's postulates are present, the basic logical postulates,
extensionality, description, infinity and choice, though not always
in Church's formulation: Dedekind's axiom stands in place of his
Axioms~7 and~8, a relational choice scheme in place of his selection
operator, and Boolean extensionality, which Church mentions but does
not postulate, is included.

The model required by the theorem can be obtained from a sufficiently
large universe type; Section~\ref{sec:harrison} makes precise how
large, and how much set theory is needed to construct it.  As an
intermediate step, the interface \isa{hol\_universe} fixes a carrier
with domains, application and abstraction satisfying the frame
conditions (full function spaces, truth values, extensionality, a
domain-respecting parameter interpretation) together with
Dedekind-infinite individuals.  Nothing is assumed about the logical
constants: the interface extends \isa{lambda\_universe} of
Section~\ref{sec:tier-one}, which supplies them.  Any interpretation
of \isa{hol\_universe}, a carrier type with domains, application and
abstraction satisfying its assumptions, yields a model as required by
Theorem~\ref{thm:confullrel}, and thereby the consistency
\artB{con\_full\_HOL\_universe}{HOL\_Universe}.

\section{A Harrison-style interface and its set-theoretic realization}
\label{sec:harrison}

\subsection{The interface}

To construct an interpretation of \isa{hol\_universe}, which still
speaks of types, domains, application and abstraction, we follow
Harrison.  His self-verification of HOL Light~\cite{Harrison06}
adds to HOL Light one axiom $I$: a type $I$ of cardinality strictly above that of the
individuals, such that every set of cardinality below $I$ has its
power set of cardinality below $I$.  In cardinal terms this is a
strong limit $\iota > \aleph_0$, for which Harrison notes the
cardinality of $\Vs_{\omega+\omega}$ as a realization.  Our rendering
is the interface \isa{harrison\_universe}, which replaces cardinal
comparison by a smallness predicate closed under encoded power
collections.  It fixes only four constants on a carrier: an injective pairing
$\mathit{pr}$, an injective encoding $\mathit{enc}$ of the \emph{small}
collections, a smallness predicate $\mathit{Sm}$, and a small
Dedekind-infinite collection $D_0$ of individuals.  These are subject to
seven assumptions: injectivity of pairing and
encoding, downward closure of smallness, closure under products, closure
under \emph{encoded power collections} (the collection of encodings of small
subcollections of a small collection is small), smallness of $D_0$, and Dedekind infinitude of
$D_0$.  It mentions neither types nor domains nor application.  Since
the interface is the development's most reusable component, we
reproduce it in full:

\begin{alltt}\footnotesize
\textbf{locale} harrison_universe =
  \textbf{fixes} pr :: "'u \(\Rightarrow\) 'u \(\Rightarrow\) 'u"    \textbf{and}  enc :: "('u \(\Rightarrow\) bool) \(\Rightarrow\) 'u"
    \textbf{and} Sm :: "('u \(\Rightarrow\) bool) \(\Rightarrow\) bool"  \textbf{and}  D0 :: "'u \(\Rightarrow\) bool"
  \textbf{assumes} pr_inj: "pr a b = pr c d \(\Longrightarrow\) a = c \(\wedge\) b = d"
    \textbf{and} enc_inj:  "\(\llbracket\)Sm X; Sm Y; enc X = enc Y\(\rrbracket\) \(\Longrightarrow\) X = Y"
    \textbf{and} Sm_sub:   "\(\llbracket\)Sm X; Y \(\le\) X\(\rrbracket\) \(\Longrightarrow\) Sm Y"
    \textbf{and} Sm_prod:  "\(\llbracket\)Sm X; Sm Y\(\rrbracket\) \(\Longrightarrow\) Sm (\(\lambda\)z. \(\exists\)a b. X a \(\wedge\) Y b \(\wedge\) z = pr a b)"
    \textbf{and} Sm_pow:   "Sm X \(\Longrightarrow\) Sm (\(\lambda\)z. \(\exists\)Y. Y \(\le\) X \(\wedge\) Sm Y \(\wedge\) z = enc Y)"
    \textbf{and} Sm_D0:    "Sm D0"
    \textbf{and} D0_inf:   "\(\exists\)f c. (\(\forall\)a. D0 a \(\longrightarrow\) D0 (f a))
                 \(\wedge\) (\(\forall\)a b. D0 a \(\longrightarrow\) D0 b \(\longrightarrow\) f a = f b \(\longrightarrow\) a = b)
                 \(\wedge\) D0 c \(\wedge\) (\(\forall\)a. D0 a \(\longrightarrow\) f a \(\ne\) c)"
\end{alltt}

\begin{theorem}[The Harrison-style interface suffices]\label{thm:harrison}
From every interpretation of \isa{harrison\_universe} an
interpretation of \isa{hol\_universe} is constructed: domains are collections of encoded functional
graphs, abstraction is encoding, application is decoding.  Consequently the
consistency of Theorem~\ref{thm:confullrel} follows from the
Harrison-style assumptions alone.
\artB{harrison\_hol\_universe}{HOL\_Universe}
\end{theorem}

In effect, the Harrison-style interface postulates one carrier type in which every
domain of the $\ii,\oo$-type hierarchy can be placed.  It is closed
under the operations that building each single type requires, pairing
and encoded power collections of \emph{small} collections, but
\emph{not} under its own full function space, which by Cantor's
theorem would be strictly larger; the carrier as a whole is
accordingly not small \artB{not\_Sm\_UNIV}{HOL\_Universe}.  This
restricted closure is a comparably weak assumption: neither
replacement nor union is demanded.  Which cardinal
hypothesis the seven assumptions correspond to exactly is left open
here.  Note that we deliberately call the interface
Harrison-\emph{style} and \emph{sufficient}: no minimality claim is made or formalized, and no
formal equivalence with Harrison's axiom $I$ is asserted; \isa{Sm\_pow}
is the counterpart of his power-set clause and \isa{Sm\_D0} that of his
clause on the individuals.  The
interface serves two purposes.  It isolates the strong-limit core of
Harrison's universe axiom $I$, the encoded power collections
of \isa{Sm\_pow}, as the assumption that carries the additional
strength of our construction.  And, mentioning neither types nor
domains nor application, it is reusable verbatim: any consistency
argument for an embedded logic that requires full function spaces
over an infinite domain of individuals can instantiate it, with Theorem~\ref{thm:harrison}
supplying the model construction.

\subsection{The set-theoretic realization}\label{sec:vmodel}

Paulson's \isa{ZFC\_in\_HOL}~\cite{ZFCinHOL} axiomatizes a
Zermelo--Fraenkel universe with choice as an Isabelle/HOL type $\Vs$.  Over it, the
domains of the intended standard model are defined by three equations:
\[
  \D_\ii = \omega, \quad
  \D_\oo = \{0, 1\}, \quad
  \D_{\sigma\Rightarrow\tau} = \prod_{x \in \D_\sigma} \D_\tau
  \;\;\text{(the full $\Vs$-function space).} \] Application is the
set-theoretic application of a function to an argument,
$\lambda$-abstraction the formation of a function's graph over the
domain, and the parameter interpretation assigns to every parameter of a type
one fixed element of that type's domain
\artB{VJv\_param\_independent}{HOL\_in\_HOL\_Deep\_Infinity\_ZFC}.
Interpreting all parameters alike is a simplification, not a
necessity: the consistency theorems below hold for any
domain-respecting parameter interpretation.

\begin{theorem}[Realization]\label{thm:vmodel}
The $\Vs$-structure satisfies the assumptions of
\isa{harrison\_universe}, with $\mathit{Sm}$ as Paulson's
smallness, encoding as set comprehension and $D_0 = \omega$, and,
directly, those of \isa{hol\_universe}, with the successor function
on $\omega$ witnessing Dedekind infinitude.
\artB{V\_harrison\_universe}{HOL\_in\_HOL\_Deep\_Infinity\_ZFC},
\artB{V\_hol\_universe}{HOL\_in\_HOL\_Deep\_Infinity\_ZFC}
\end{theorem}

\begin{theorem}[Consistency of full classical HOL]\label{thm:confull}
$\;\mathrm{con}\,\bigl(\{\DInf\} \cup
  \{\ACrel_{\sigma\tau} \mid \sigma, \tau\}\bigr)$, at every signature.
\artB{con\_full\_HOL}{HOL\_in\_HOL\_Deep\_Infinity\_ZFC}
\end{theorem}

Theorem~\ref{thm:confull} is Theorem~\ref{thm:harrison} applied to
the instance of Theorem~\ref{thm:vmodel}; it depends on the
$\Vs$-structure only through the seven assumptions of the
Harrison-style interface.  In particular $\mathrm{con}\,\{\DInf\}$
holds at every signature
\artB{con\_DInf}{HOL\_in\_HOL\_Deep\_Infinity\_ZFC}, and this is
exactly the premise of the conditional model existence of
Section~\ref{sec:completeness-vs-standard}.  Combining the
two yields a \emph{countable} Henkin model of $\DInf$ whose
construction never mentions $\Vs$; the set-theoretic universe serves
only to discharge the premise.

\begin{remark}[Readings at lower strength]\label{rem:foundations}
Theorem~\ref{thm:confull} is a theorem of Isabelle/HOL extended by
\isa{ZFC\_in\_HOL}; it is therefore a consistency result
\emph{relative} to that theory.  Two readings at lower strength are
metatheoretic, made precise in the documentation of $\mathcal{A}_2$
and not themselves formalized theorems.  First, choice can be
dropped: G\"odel's constructible universe interprets ZFC in ZF, so
the result holds relative to ZF.  Second, the construction uses no
replacement, so it can be carried out in Zermelo set theory with
choice; whether choice can be dropped there as well is
delicate~\cite{Mathias01}.  The rank analysis of the next section,
which places every domain in $\Vs_{\omega+\omega}$, is the formal
basis of both readings.
\end{remark}

\section{Rank analysis}\label{sec:rank}

We now bound how much of the set-theoretic universe the model of
the previous section occupies.  Recall that the von Neumann universe
is built in stages, $\Vs_0 = \emptyset$,
$\Vs_{\alpha+1} = \mathcal{P}(\Vs_\alpha)$ and
$\Vs_\lambda = \bigcup_{\alpha<\lambda} \Vs_\alpha$ at limits
$\lambda$, and that the rank of a set is the first stage at which it
appears as a subset.  $\Vs_\omega$ is the union of the finite stages, and
$\Vs_{\omega+\omega}$ the union of the stages $\Vs_{\omega+n}$
obtained from $\Vs_\omega$ by finitely many power-set steps; these
are the sets Zermelo set theory provides.

\begin{theorem}[Rank bound, type by type]\label{thm:rank}
In the concrete model, $\mathrm{rank}(\D_\ii) = \omega$, and for every
type $\sigma$:
\[
  \mathrm{rank}(\D_\sigma) < \omega + \omega,
  \qquad\text{hence}\qquad
  \D_\sigma \in \Vs_{\omega+\omega}.
\]
\artB{rank\_D\_iota}{HOL\_in\_HOL\_Deep\_Infinity\_ZFC},
\artB{rank\_D\_lt}{HOL\_in\_HOL\_Deep\_Infinity\_ZFC},
\artB{D\_in\_Vset\_ww}{HOL\_in\_HOL\_Deep\_Infinity\_ZFC}
\end{theorem}

Each function-space layer is contained in a triple power set of the
union of its argument and result domains, so it raises the rank by at
most three power-set steps above those domains; the limit ordinal
$\omega + \omega$ absorbs these finite increments.  The statement is one theorem, quantified over the object-level type $\sigma$, a value of the
datatype of types.  The set of all domains,
$\{\D_\sigma \mid \sigma \text{ a type}\}$ taken as one object, has
rank $\omega + \omega$, since the domains' ranks grow unboundedly
below $\omega+\omega$, and is therefore \emph{not} an element of
$\Vs_{\omega+\omega}$ (not formalized).  The bound substantiates the readings of the previous
section: a single derivation mentions finitely many types, so its
refutation uses only sets of bounded rank, finitely many power-set
steps above $\omega$.

\section{Lessons from the formalization}\label{sec:lessons}

\subsection*{Two tiers, one boundary}
Soundness alone shows where the first tier ends: finite models satisfy
the first-order sentences and refute $\DInf$, so finite-model
arguments certify the sentences but not the axiom.  Soundness alone
does not show that the two senses of infinity differ: at standard
models they coincide (Proposition~\ref{prop:dinfsat}), and separating
them needs a model of all the sentences that refutes $\DInf$.  Such a
model exists only in Henkin semantics, where function domains may be
sparse, and it is the completeness proof that builds it
(Theorem~\ref{thm:henkinsep}, Remark~\ref{rem:semantic-gap}).

\subsection*{Eigenparameters and finitary consequence}
Over a fixed signature, an infinite context may use every parameter
there is.  The quantifier rule $\NK(\Pi I)$ then cannot be applied,
because its eigenparameter, the parameter standing for an arbitrary
element of its type, must be fresh for the context.\footnote{A
  by-product of formulating $\NK(\Pi I)$ with eigenparameters; a
  formulation over free variables would avoid the restriction and is
  planned for the next version of $\mathcal{A}_1$.}  Andrews-style
finitary consequence $\fder$, where some finite part of the context
derives the conclusion, avoids the restriction as well: a finite
context leaves fresh parameters over an infinite signature, so no
purity condition is needed.

\subsection*{Abstract interfaces before concrete models}
Interposing \isa{hol\_universe} and \isa{harrison\_universe} between
the consistency theorem and $\Vs$ separates what the proof needs from
where it is realized: Theorem~\ref{thm:harrison} depends on seven
assumptions about four constants, and $\Vs$ satisfies these
assumptions but is itself not part of the proof.  The interfaces also
make dependencies visible in what they do not mention: two properties
stem from the meta-logic rather than from the interfaces.  The
object-level description operator is interpreted by Isabelle's
definite description \isa{THE}, and the relational choice scheme holds
in every standard model because the meta-logic's own choice supplies
the witness (Section~\ref{sec:choice}).

Another design decision is the locally nameless syntax of
$\mathcal{A}_1$: no proof in the
development ever renames a bound variable, in exchange for a thin
layer of named binders in which the axiom $\DInf$ of
Section~\ref{sec:boundary} is written.

\section{Related work}\label{sec:related}

This paper describes two AFP entries: $\mathcal{A}_1$~\cite{HOLinHOLDeep},
which formalizes BKK's syntax, Henkin semantics and
calculus~\cite{BBK04} with Andrews-style description~\cite{Andrews72}
and was extended for this paper, and its companion
$\mathcal{A}_2$~\cite{InfinityZFC}, which adds the second tier.  Closest in method, a deep embedding of a Church-style type theory
with machine-checked metatheory, is the work on Andrews'
$\mathcal{Q}_0$.
D\'iaz~\cite{Q0MetatheoryAFP} machine-checks its soundness and, by
Andrews' one-element standard model (Theorem~5403
of~\cite{Andrews02}), its consistency; Schlichtkrull's
entry~\cite{Q0SoundnessAFP}, developed independently in the same
month, formalizes soundness as well.  Building on D\'iaz's entry,
Boserup, Huerta y Munive and Schlichtkrull add Henkin
completeness~\cite{Q0CompletenessAFP} by an abstract consistency
property, using a locale framework for such
properties~\cite{FromSchlichtkrull25}; see also Koch and Kirst's
analysis of second-order logic in Coq~\cite{KochKirst22}.  The $\mathcal{Q}_0$ line develops its
semantics over \isa{ZFC\_in\_HOL} throughout, the one-element model
included; the present embedding admits set theory only on the second
tier, so that tier one, the same classical argument for $\NK$, is a
theorem of plain Isabelle/HOL with an audited footprint.
Closest in question is Harrison's self-verification of HOL
Light~\cite{Harrison06}: machine-checked consistency of HOL with
infinity and choice is not new, nor is the division into two tiers.
Harrison proves soundness and
consistency of HOL Light's core twice, without its axiom of infinity
inside HOL Light and with it inside HOL Light plus a strong-limit
universe axiom, and already notes $\Vs_{\omega+\omega}$ as a realization.  Its
successors by Kumar, Arthan, Myreen and
Owens~\cite{Kumar16} and the verified implementation
Candle~\cite{Candle} extend the soundness result to HOL's definitional
mechanisms, taking the existence of a set-theoretic universe as an
assumption of the soundness theorem.  All of these work with standard
semantics and prove soundness; none has a completeness theorem.
Harrison's object logic has polymorphism and the choice operator
$\varepsilon$, ours has primitive equality and typed description, so
neither consistency result subsumes the other.  Beyond completeness itself, the following are, to our knowledge,
new.  On the Henkin side: conditional model existence on the first
tier, in hypothesis-relative form (Sections~\ref{sec:embedding}
and~\ref{sec:completeness-vs-standard}); the characterization of
satisfaction of $\DInf$ in Henkin models and the measurement of the
boundary between the two senses of infinity by $\NK$-derivations and
a Henkin countermodel (Section~\ref{sec:infinity}); the countable
Henkin model of $\DInf$ for which the universe discharges only the
premise $\mathrm{con}\,\{\DInf\}$ (Section~\ref{sec:vmodel}).  On the
standard side: the per-type rank bound as a theorem
(Section~\ref{sec:rank}).  Refinements of Harrison's two tiers are, on
the first, a possibly finite signature, no choice operator in the
definition of the model and a mechanical footprint audit (Section~\ref{sec:tier-one}),
and, on the second, the factorization of the argument through two
reusable locales, one of which isolates a version of Harrison's axiom
(Sections~\ref{sec:universe} and~\ref{sec:harrison}).

The object logic is also the base of dependently typed HOL
(DHOL)~\cite{DHOLCADE,DHOL}, which translates into HOL by erasing the
dependencies, soundly and completely for provability~\cite{DHOLCADE},
so that the consistency results of this paper transfer to it.  Rabe's
lax models for DHOL~\cite{RabeDHOL}, a Henkin-style generalization of
standard models, come with soundness and completeness theorems, so
the boundary drawn here can be drawn for DHOL as well (see the
conclusion).  For Lean's theory, by contrast, Carneiro gives a
set-theoretic model and thereby consistency relative to ZFC with
countably many inaccessible cardinals~\cite{Carneiro19}, each universe level modelling the theory with one level fewer; the
metatheory of its definitional equality is still being
settled~\cite{Carneiro24}.  Every level contains the natural numbers,
so it contains no infinity-free tier with finite models, and there
is no completeness theorem that ties provability to a class of
models.

Nipkow and Ro{\ss}kopf formalize Isabelle's own
metalogic with a verified proof checker~\cite{NipkowRosskopf21}, the
same programme for the logical framework beneath Isabelle/HOL.

\section{Conclusion}\label{sec:conclusion}

What a copy of higher-order logic can certify about itself has two
tiers, measured by the strength the certifying copy needs, and both
correspond to the two results of Harrison's self-verification of HOL
Light.  On the first tier, plain Isabelle/HOL certifies Church's type
theory without an axiom of infinity, by a standard model with one
individual.  On the second, a Harrison-style universe axiom realized
in \isa{ZFC\_in\_HOL} suffices for full classical higher-order logic
with infinity and relational choice, with the domain of every type of
rank below $\omega+\omega$.  What the present formalization adds is
Henkin completeness, and with it a semantics that matches the calculus
exactly and, on the first tier, compactness and model existence
conditional on consistency alone.  Even the countable Henkin model of
the Dedekind axiom of infinity is therefore built in plain
Isabelle/HOL.
The boundary between the tiers runs between sentences that have
finite models and sentences that have none; the two senses of
infinity are separated by a countable Henkin model that satisfies
every sentence ``there are at least $n$ individuals'' and refutes the
axiom.  Every theorem, proposition and corollary stated above is a
machine-checked theorem of the two AFP entries $\mathcal{A}_1$ and
$\mathcal{A}_2$; the metatheoretic readings are marked as not
formalized.

Three directions remain open.  Whether the Harrison-style interface is
minimal, in particular whether the encoded-power assumption
\isa{Sm\_pow} is necessary for full function spaces, is not settled
here.  The equiconsistency of Isabelle/HOL's monomorphic core with
$\NK$ plus infinity and choice is a candidate for formalization in
its own right; on that basis, step~(4) of the argument in
Section~\ref{sec:boundary}, the appeal to G\"odel's second
incompleteness theorem, could be formalized as well.  It would require
interpreting enough arithmetic inside $\NK + \DInf$, with Paulson's
mechanized incompleteness theorems~\cite{Paulson15} as the model to
follow; the natural first step, Dedekind's own theorem that the axiom yields a model of the
Peano axioms, would itself be a worthwhile derivation inside $\NK$.

Finally, the two interfaces of Sections~\ref{sec:universe}
and~\ref{sec:harrison} invite reuse beyond this paper.  The natural
first candidate is dependently typed HOL (Section~\ref{sec:related}):
connecting a formalization of Rabe's DHOL semantics~\cite{RabeDHOL} to
the interfaces of this paper would test them on a logic properly
beyond Church's.

\subsection*{Artifacts}
Both artifacts are entries of the Archive of Formal Proofs, released
with Isabelle2026 in October 2026, and both build with Isabelle2026.
$\mathcal{A}_1$ is \isa{HOL\_in\_HOL\_Deep}~\cite{HOLinHOLDeep}
(\url{https://isa-afp.org/entries/HOL_in_HOL_Deep.html}), first
published in August 2026 and extended in September 2026 for this
paper.  $\mathcal{A}_2$ is its companion
\isa{HOL\_in\_HOL\_Deep\_Infinity\_ZFC}~\cite{InfinityZFC}
(\url{https://isa-afp.org/entries/HOL_in_HOL_Deep_Infinity_ZFC.html}),
which also uses the AFP entry \isa{ZFC\_in\_HOL}.  The versions this
paper describes extend the released entries by a few facts cited
above (\isa{DInf\_sat\_iff\_Dedekind}, \isa{not\_Sm\_UNIV},
\isa{countable\_henkin\_model\_DInf\_of\_con}, and the direct
finite-model proofs of Section~\ref{sec:infinity}, which drop the
infinite-signature hypothesis from Theorems~\ref{thm:distinctcon}
and~\ref{thm:distinctnotdinf}) and by the theory
\isa{Footprint\_Deps} of Section~\ref{sec:tier-one}.  The complete sources in the versions
described here accompany the arXiv version as ancillary files, and
every artifact pointer in this paper resolves in them.  The same
extensions are part of the Archive's development version since
November 2026 (\url{https://devel.isa-afp.org/entries/HOL_in_HOL_Deep.html},
\url{https://devel.isa-afp.org/entries/HOL_in_HOL_Deep_Infinity_ZFC.html})
and will be in its next release.  The theory \isa{Main\_Results} of
$\mathcal{A}_1$ restates the principal results with minimal hypotheses
and without internal notation and is the intended entry point;
$\mathcal{A}_2$ is read in order, \isa{HOL\_Universe} (the two
interfaces) followed by the $\Vs$-instantiation.
Table~\ref{tab:index} indexes the principal results against their
formal counterparts.  $\mathcal{A}_1$ comprises eleven theories, about
7{,}500 lines and 462 named facts, and imports only the Isabelle
distribution; $\mathcal{A}_2$ comprises two theories, about 950 lines
and 69 named facts, and imports $\mathcal{A}_1$ and
\isa{ZFC\_in\_HOL}.  Every proof in both artifacts, including the 23 in $\mathcal{A}_1$
that call the \isa{smt} method, is replayed through Isabelle's
kernel.  Software Heritage identifiers
follow with the final version, as the journal's artifact policy asks.

\begin{table}[t]
\caption{Index of the principal results; \isa{Infinity\_ZFC} abbreviates
  the theory \isa{HOL\_in\_HOL\_Deep\_Infinity\_ZFC} of $\mathcal{A}_2$.}
\label{tab:index}
\scriptsize
\setlength{\tabcolsep}{3pt}
\begin{tabular}{@{}l>{\raggedright\arraybackslash}p{0.29\linewidth}>{\raggedright\arraybackslash}p{0.35\linewidth}>{\raggedright\arraybackslash}p{0.22\linewidth}@{}}
\toprule
result & & Isabelle fact & theory \\
\midrule
Thm~\ref{thm:soundcomplete} & soundness and completeness &
  \isa{sound\_and\_complete\_hyps} & \isa{Main\_Results} ($\mathcal{A}_1$) \\
Cor~\ref{cor:compact} & semantic compactness &
  \isa{NK\_consequence\_compact} & \isa{Main\_Results} ($\mathcal{A}_1$) \\
Thm~\ref{thm:consistency} & consistency without infinity &
  \isa{consistency} & \isa{Main\_Results} ($\mathcal{A}_1$) \\
Prop~\ref{prop:dinfref} & finite refutation &
  \isa{DInf\_refuted\_finite} & \isa{NK\_Infinity} ($\mathcal{A}_1$) \\
Thm~\ref{thm:dinfdistinct} & axiom $\vdash$ every $\Exn{n}$ &
  \isa{axiom\_derives\_\allowbreak every\_\allowbreak finite\_cardinality} & \isa{Main\_Results} ($\mathcal{A}_1$) \\
Prop~\ref{prop:dinfsat} & satisfaction of $\DInf$ &
  \isa{DInf\_sat\_iff}, \isa{DInf\_sat\_iff\_Dedekind} & \isa{NK\_Infinity} ($\mathcal{A}_1$) \\
Thm~\ref{thm:distinctcon} & the sentences are consistent &
  \isa{finite\_cardinalities\_\allowbreak consistent} & \isa{Main\_Results} ($\mathcal{A}_1$) \\
Thm~\ref{thm:distinctnotdinf} & $\{\Exn{n}\}_n \nVdash$ axiom &
  \isa{finite\_cardinalities\_\allowbreak do\_not\_\allowbreak derive\_axiom} & \isa{Main\_Results} ($\mathcal{A}_1$) \\
Thm~\ref{thm:henkinsep} & Henkin countermodel &
  \isa{henkin\_scheme\_refutes\_DInf} & \isa{NK\_Infinity} ($\mathcal{A}_1$) \\
Thm~\ref{thm:completeness-dinf} & completeness of the extension &
  \isa{NK\_completeness\_DInf} & \isa{HOL\_Universe} ($\mathcal{A}_2$) \\
Thm~\ref{thm:confullrel} & abstract consistency &
  \isa{con\_full\_HOL\_rel} & \isa{HOL\_Universe} ($\mathcal{A}_2$) \\
Prop~\ref{prop:acrelsat} & choice at standard models &
  \isa{ACrel\_sat} & \isa{HOL\_Universe} ($\mathcal{A}_2$) \\
Thm~\ref{thm:harrison} & Harrison-style interface suffices &
  \isa{harrison\_hol\_universe} & \isa{HOL\_Universe} ($\mathcal{A}_2$) \\
Thm~\ref{thm:vmodel} & $\Vs$-realization &
  \isa{V\_harrison\_universe} & \isa{Infinity\_ZFC} ($\mathcal{A}_2$) \\
Thm~\ref{thm:confull} & consistency of full HOL &
  \isa{con\_full\_HOL} & \isa{Infinity\_ZFC} ($\mathcal{A}_2$) \\
Thm~\ref{thm:rank} & type-by-type rank bound &
  \isa{rank\_D\_lt}, \isa{D\_in\_Vset\_ww} & \isa{Infinity\_ZFC} ($\mathcal{A}_2$) \\
\bottomrule
\end{tabular}
\end{table}

\subsection*{Disclosure}
The argument, structure and final text of this manuscript are the
authors' own.  A generative AI assistant (Anthropic's Claude) was
used, under continuous human direction, for three kinds of tasks:
drafting Isabelle proof scripts and their documentation, all of which
the authors reviewed, revised and simplified; editorial revision of
passages of this manuscript according to the authors' annotations;
and checking the correspondence between the informal statements and
the Isabelle statements, as well as bibliographic data.  All proofs
are machine-checked by Isabelle, and the authors take full
responsibility for the content.  This work received no external
funding, including no subsidized access to tools or software.

\bibliographystyle{amsplain}
\bibliography{afm_paper}

\end{document}